\documentclass[%
 reprint,
 amsmath,amssymb,
 aps,
prd,
]{revtex4-2}

\usepackage{graphicx}
\usepackage{dcolumn}
\usepackage{color}
\usepackage{bm}

\def\nn{\nonumber}

\def\dA{{\rm d}A}

\def\R2{{\mathcal R}^2}

\def\UU{\displaystyle \cup }
\def\btc{\begin{tcolorbox}}
\def\etc{\end{tcolorbox}}
\def\be{\begin{equation}}
\def\ee{\end{equation}}
\def\barr{\begin{array}{lr}}
\def\earr{\end{array}}
\def\bea{\begin{eqnarray}}
\def\eea{\end{eqnarray}}
\def\la{\langle}
\def\ra{\rangle}
\def\vx{{\mathbf x}}

\def\nn{\nonumber}

\def\a{\alpha}

\def\s{\sigma}

\def\bs{b_{\rm sum}}

\begin{document}
\title{On the statistical nature of Betti numbers and Euler characteristic of smooth random fields}

\author{Pravabati Chingangbam$^{1,2}$} 
\email{prava@iiap.res.in}
\affiliation{$^1$ Indian Institute of Astrophysics, Koramangala II Block,      
  Bengaluru  560 034, India\\
  $^2$ School of Physics, Korea Institute for Advanced Study, 85 Hoegiro, Dongdaemun-gu, Seoul 02455,  Republic of  Korea
  }
\begin{abstract}
We represent excursion sets of smooth random fields 
as  unions of a topological basis consisting of a sequence of simply and multiply connected compact subsets of the underlying manifold. The associated coefficients, which are non-negative discrete random variables, reflect the randomness of the field.  
Betti numbers of the excursion sets can be expressed as summations over the coefficients, and the Euler characteristic and the sum of Betti numbers can also be expressed as their (alternating) sum. This enables understanding their statistical properties as sums (or differences) of discrete random variables. 
We examine the conditions under which each topological statistic can be asymptotically Gaussian as the size of the manifold and the resolution increase.
The coefficients of the basis elements are then modeled as Binomial variables, and the statistical natures of Betti numbers, Euler character and sum of Betti numbers follow from this fundamental property. We test the validity of the modeling using numerical calculations, and identify threshold regimes where the topological statistics can be approximated as Gaussian variables.  
The new representation of excursion sets thus maps the properties of topological statistics to combinatorial structures, thereby providing mathematical clarity on their use for physical inference, particularly in cosmology.
\end{abstract}

\maketitle

\section{Introduction}
\label{sec:s1}

Random fields are used to describe a wide variety of phenomena in nature. In cosmology, observed spatial data are modeled as random fields, and physical inference is drawn using statistics constructed from the fields. 
Examples of such statistics are the power spectrum, bispectrum, higher $n$-point functions; and topological and geometrical statistics~\cite{Adler:1981,Adler:2007}. 

This note is concerned with topological statistics, namely, Betti numbers and Euler characteristic that quantify the topological information  of excursion sets of random fields. 
Early applications of topological statistics to cosmology focused on using the genus to characterize the large scale structure of the universe~\cite{Gott:1986APJ,Hamilton:1986}. The genus is equivalent but not equal to the Euler characteristics in terms of structural information.  Minkowski functionals (MFs), one of which is the Euler characteristic, were introduced to study cosmological fields~\cite{Mecke:1994,Schmalzing:1997,Gott:1990,Schmalzing:1998,Ducout:2013}, thereby exploiting geometrical properties of excursion sets, in addition to their topology. 
Subsequently, Betti numbers, whose alternating sum gives the Euler characteristic, were also added to the suite of statistics. They have been studied numerically for excursion sets of Gaussian random fields in~\cite{Park:2013,Pranav:2017,Pranav:2019,Feldbrugge:2019}. They are studied under the ambit of topological data analysis and persistent homology~\cite{Edelsbrunner:2022computational,Sousbie:2011persistent,Weygaert:van2011,Weygaert:2013,Bermejo:2024}.  In~\cite{Chingangbam:2012} they were used to probe non-Gaussianity of the cosmic microwave background (CMB), and it was shown that they can provide more information compared to Euler characteristic (see also~\cite{Cole:2018,Feldbrugge:2019,Biagetti:2020}). 
More recent applications of geometrical and topological statistics include cosmological parameter inference~\cite{Appleby:2020genus,Appleby:2022,Liu:2022,Liu:2023}, probing statistical properties of the CMB~\cite{Pranav:2022,Pranav:2023,Bashir:2025}, understanding Galactic foreground components~\cite{Chingangbam:2013,Rana:2018oft,Rahman:2021azv,Rahman:2022,Martire:2023}, constraining the physics of the epoch of reionization~\cite{Kapahtia:2017qrg,Bag:2018,Kapahtia:2019ksk,Akanksha:2021,Giri:2021,Wilding:2021,Elbers:2023}, understanding interstellar magnetic fields \cite{Shukurov:2018}, 
and inter-Galactic medium~\cite{Chepurnov:2008}, probing Galactic properties~\cite{Spina:2021} and Cosmic web~\cite{Tymchyshyn:2023czh}, topology of the universe~\cite{Aurich:2024}, and validation of machine learning  methods~\cite{Grewal:2024}.

Minkowski functionals and Betti numbers of excursion sets  are random variables, owing to the randomness of the field. In principle, their statistical properties, namely,  ensemble expectations, variances and higher moments, and probability density functions, must be derived from the statistical nature of the field. 
Analytic formulas for these statistical properties, when known,  are invaluable because they serve to validate and guide numerical computations in practical applications. In this regard, considerable theoretical progress has been made in calculating analytic formulas for ensemble expectations of MFs for  fields that are Gaussian~\cite{Adler:1981,Tomita:1986}, lognormal~\cite{Coles:1988,Coles:1991}, Rayleigh~\cite{Naselski:1998} and mildly non-Gaussian~\cite{Matsubara:2003,Matsubara:2011,Matsubara:2020,Gay:2012,Chingangbam:2024}.  
In comparison, Betti numbers are much harder to study analytically and obtaining analytic formulas for their ensemble expectations, even for Gaussian random fields, remains an open problem. 

Efforts to understand the statistical behavior of Minkowski functionals and Betti numbers beyond their ensemble averages remain largely absent in the cosmology literature. It is often assumed - without rigorous analytic justification - that these quantities follow Gaussian statistics. Typically, their standard deviations, derived from simulations, are employed to assess the statistical significance of measurements obtained from observational data. However, the increasing availability and precision of cosmological data have expanded the range of applications that demand a more accurate characterization of their statistical properties. Consequently, broadening our understanding in this area has become imperative. This growing interest is also reflected in recent advances within the mathematics community~\cite{Owada:2020,Owada:2023}.

Toward this goal, in this paper we  
investigate the statistical nature of Betti numbers, Euler characteristic and the sum of Betti numbers for random fields on finite sized manifolds. We begin by representing excursion sets  
  as unions of a topological basis consisting of a sequence of multiply connected compact subsets of the manifold. 
  This representation allows a clean separation of the probabilistic nature of the field encoded in the coefficients, from the topological information  of the basis. 
By expressing  Betti numbers, Euler characteristic and the sum of Betti numbers in terms of the basis coefficients, we are able to formulate the question of their statistical properties in a mathematically precise manner. 
The paper is organized as follows. Sec.~\ref{sec:s2} gives a brief description of smooth random fields and physically relevant scales. The definitions of Betti numbers and Euler characteristic are given in Sec.~\ref{sec:s3}. Sec.~\ref{sec:s4} introduces the new topological basis for excursion sets, and discusses the statistics of  Betti numbers, Euler characteristic and sum of Betti numbers in terms of coefficients of the basis in two dimensions. The possibility of their  asymptotic Gaussian behavior is discussed. Then, the coefficients are argued to be  Binomial variables and the implications of this on the statistical nature of Betti numbers, Euler characteristic and sum of Betti numbers are investigated.  Generalization to three dimensions is carried out in Sec.~\ref{sec:s5}. Lastly, the results are summarized along with a discussion of the broader implications in Sec.~\ref{sec:s6}.

\section{Smooth random fields and excursion sets}
\label{sec:s2}

This section summarizes the basic definitions of random fields and excursion sets, and the relevant length scales  that we will require in subsequent sections.  For details we refer the reader to~\cite{Adler:2007,Yaglom:1987,BBKS}. 

Let $M$ denote a smooth $d$-dimensional manifold.  A function $f(\vx)$ is said to be a random field on $M$ if at each $\vx\in M$, $f$ is a random variable, and the joint probability density function (PDF), ${\cal P}(f(\vx_1),f(\vx_2),...,f(\vx_k))$, of $f$ at $k$ arbitrary points on $M$ is well defined. The covariance of $f$  between any two points $\vx$ and $\vx'$ gives the covariance function $\xi(\vx,\vx')$.  The auto-covariance $\xi(\vx,\vx)=\s^2_0(\vx)$ gives the variance of $f$ at $\vx$.  Higher order correlation functions can be likewise defined. The field is said to be {\em smooth} if its realizations are differentiable to all orders. The derivatives of the field are also random fields, and the field and its derivatives then form a set of multivariate random variables at each point on $M$. $f(\vx)$ is said to be  Gaussian if  ${\cal P}$ has Gaussian form. If $f$ is Gaussian, then its derivatives are also Gaussian fields.

The field is said to be (strictly) homogeneous if  ${\cal P}$  is invariant under group actions that take point  $\vx$  to $\vx + \mathbf{a}$ on $M$. This condition implies that the covariance and all higher order correlation functions depend only on the vector $\vx-\vx'$. 
(The converse need not be true.)  
Homogeneity is a necessary condition for the field to be ergodic~\cite{Adler:1981}. In cosmology, where we observe only a single realization of a random field, ergodicity is implicitly assumed when ensemble averages are replaced by volume averages to infer physical quantities from data. Further, the field is said to be isotropic if the covariance and all higher order correlation functions depend only on the distance $|\vx-\vx'|$. 

{\em Definition of excursion set}: Let $\nu\sigma_0$ denote level values of a homogeneous field $f$, and let  $\sigma_0$ denote its standard deviation. $\nu$ is usually referred to as `threshold'. The set of all points on $M$ where $f \ge\nu\sigma$, which we denote by $\mathcal{Q}^{(\nu)}$, is called the {\em excursion set} indexed by $\nu$. By varying $\nu$ we obtain a one-parameter family of excursion sets. 

{\em Spectral parameters}: 
Let the variances of the field and its $n$-th derivatives be denoted by $\s_n^2$. 
Successive ratios,   $\s_n / \s_{n+1}$, 
give a hierarchy of length scales associated with the field. Of these, the first one is $r_c=\s_0/\s_1$, and it represents the typical size of spatial fluctuations.

Let us now focus on an isotropic field on $M=R^3$. The spectral parameters are related to the power spectrum of the field, $P(k)$, by the relation~\cite{BBKS}, 
 \be
 \s_n^2=\int_0^{\infty} dk \frac{k^2}{2\pi^2} k^{2n}P(k).
 \label{eq:sn}
\ee
If the field is smoothed with a Gaussian kernel with smoothing scale $R_s$, then $\s_n^2$ has the form,
 \be
 \s_n^2=\int_0^{\infty} dk \frac{k^2}{2\pi^2} k^{2n}P(k) W^2(kR_s),
 \label{eq:sns}
\ee
where $W(kR_s)$ is the Fourier transform of the Gaussian kernel. Fluctuations at scales smaller than $R_s$ get washed out due to the smoothing. 
From Eqs. \ref{eq:sn} and \ref{eq:sns} we get
\be r_c^2=\frac{\int dk k^2 P(k) W^2(kR_s)}{ \int dk k^4 P(k) W^2(kR_s) }.
\ee

For cosmological data there are two length scales introduced by the observation process. The first is the  size, $L$, of  $M$ which represents the finite field of view (usually a two or three dimensional compact manifold). The second is the resolution of the data which introduces a minimum length scale. We roughly identify this scale with $R_s$. We work with the condition that $L\gg R_s$.  Information of both $L$ and $R_s$ are folded in $r_c(L,R_s)$ via low and high mode cut-offs in the integrals in Eqs.~(\ref{eq:sn}) and (\ref{eq:sns}).

For subsequent sections it is useful to discuss how $r_c$ scales relative to $L$ \footnotemark[1]\footnotetext[1]{Here $L$ increases due to increase of the field of view. We are not considering scaling transformations of the manifold.} and $R_s$ which are observer dependent, based on the form of $P(k)$. 
Keeping the discussion general,  we  broadly distinguish two forms of $P(k)$:
\begin{enumerate}
\item $P(k)$  does not contain any intrinsic physical cut-off scale other than $L$ and $R_s$. An example is a power law form with spectral index $\a$. 
  This case is of interest in cosmology because cosmological fields are assumed to have evolved from primordial density fluctuations with nearly scale invariant power spectrum  within the inflationary paradigm. 
In this case, for fixed $L$, we have $r_c\to 0$ as  $R_s\to 0$. On the other hand, for fixed $R_s$, we have $r_c\to \infty$ as $L\to \infty$.
\item $P(k)$ has a form that contains intrinsic small scale cut-off at some $k=k_s > R_s^{-1}$, or  large scale cut-off at some $k=k_L< L^{-1}$, or both.  Such examples arise when some physical process erases small and/or lager scale fluctuations.   In this case we have  $r_c=r_c(k_L^{-1}, k_p^{-1})$, and it is independent of $L, R_s$.  
\end{enumerate}

Next, let us define a `packing fraction' parameter $q$ as,
\be
 q\equiv \left(\frac{L}{r_c}\right)^d.
  \label{eq:q}
  \ee
  $q$ roughly gives the number of `structures', of size $r_c$, that can fit in a region of size $L$ in $d$-dimensional space.  The meaning of structure will be made more precise in subsequent sections, and here we can assume it means a circle in two dimensions, or a sphere in three dimensions. We are interested in the question of whether $q$ can approach infinity in the limits of $L\to\infty$ and $R_s\to 0$. The scaling behaviour of $q$ is summarized in table \ref{tab:t1}.

\renewcommand{\arraystretch}{1.2}
\begin{table}[ht]
  \begin{tabular}{|p{3.cm}|l|l|l|} \hline
Limits  & \ $P(k)$ & $\ r_c$ &\ $q$  \\ \hline 
  $L\to \infty$, $R_s$ fixed  & 
\begin{tabular}{c}
  Type 1 \\
  Type 2
\end{tabular}       &
\begin{tabular}{l}
  $r_c \propto L$\\
   $r_c $ constant 
\end{tabular}    &
\begin{tabular}{l}
  $q $ constant\\
   $q \to \infty$ 
\end{tabular}    \\
\hline
  $L$ constant, $R_s\to 0$  & 
\begin{tabular}{c}
  Type 1 \\
  Type 2
\end{tabular}       &
\begin{tabular}{l}
  $r_c \propto R_s $\\
   $r_c $ constant
\end{tabular}    &
\begin{tabular}{l}
  $q \to \infty$\\
   $q $ constant
\end{tabular}    \\
\hline
\end{tabular}
\caption{Scaling of $q$ in the limits $L\to \infty$ and $R_s\to 0$.}
\label{tab:t1}
\end{table}
Note that although Eqs.~(\ref{eq:sn}) and (\ref{eq:sns}) are valid for $M=R^3$ the argument for the asymptotic behaviour of $r_c$ can be extended to general $M$. 
In the next section we will use $q$ to set limits on the number of structures, which in turn will determine how (and whether) Betti numbers, Euler characteristic and the sum of Betti numbers approach Gaussian nature in the limit of $L\to\infty$ and $R_s\to 0$. 

\section{Definitions of Betti numbers, Euler characteristic and sum of Betti numbers}
\label{sec:s3}

Betti numbers are topological invariants that can be used to distinguish topological spaces based on the notion of connectivity (which inform whether the space is a single piece or consists of disjoint pieces) and holes.  Formally, they are defined as the ranks of the homology groups $H_k$ of a topological space. The $k$th Betti number is the rank of the $k$th homology group~\cite{Pontryagin:1952,Munkres:1984}. 
Intuitively, $k$ gives the number of $k$-dimensional `holes' in the space. 
 
In this paper we are concerned with Betti numbers of two and three dimensional topological spaces. In what follows, we discuss some examples of simple topological spaces that are relevant for subsequent sections.

\begin{figure}[ht]
  \centering $Q_0$  \hskip 2.5cm $Q_1$ \hskip 2.2cm $Q_2$ \\
  \centering\includegraphics[height=1.4cm,width=8.cm]{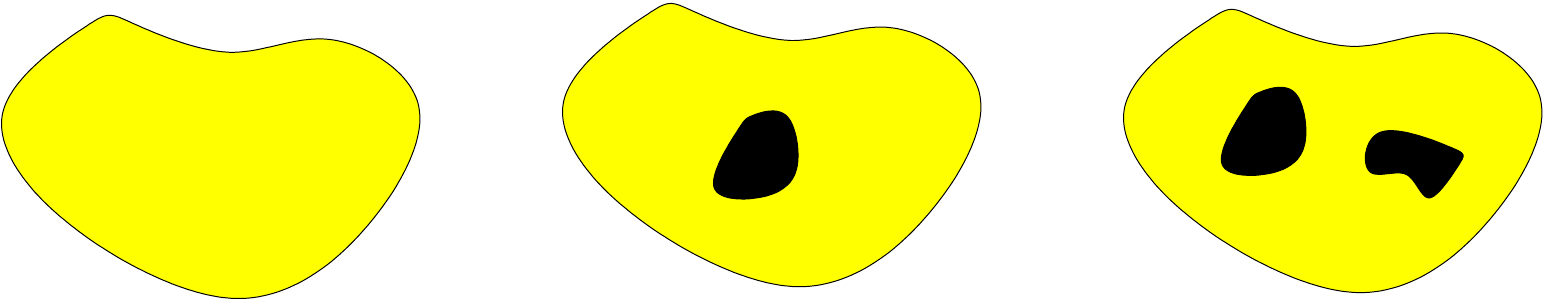}\\
\caption{Examples of basic topological spaces in two dimensions - a solid disk (left), an annular disk containing one hole shown in black (middle), and a disk containing two holes (right).}
\label{fig:2dspaces}
\end{figure}

{\em Betti numbers for two dimensions}: Let us consider two dimensional topological spaces that have orientable boundaries. There are two Betti numbers for such spaces, namely, $b_0$ which counts the number of connected components (or disjoint pieces), and $b_1$ which counts the number of one dimensional holes. The Euler characteristic, $\chi$, is given by the difference of the Betti numbers, as,
\be\chi=b_0-b_1.\ee
We also introduce the sum of Betti numbers
\be
b_{\rm sum}= b_0+b_1. 
\ee

Fig.~\ref{fig:2dspaces} shows examples of a disk denoted by $Q_0$ (left), an annular disk containing one hole denoted by $Q_1$ (middle), and a disk containing 2 holes denoted by $Q_2$  (right).  $Q_0$ has one positively (by convention) oriented boundary enclosing it. It has $b_0=1, \,b_1=0,\, \chi=1,\,\bs=1$. $Q_1$ has one positively oriented boundary and one negatively oriented one enclosing the hole. It has $b_0=1, \,b_1=1,\, \chi=0,\,\bs=2$. $Q_2$ has one positively oriented boundary and two negatively oriented ones enclosing the two holes. It has $b_0=1, \,b_1=2,\, \chi=-1,\,\bs=3$.  For our purpose here, topological spaces that have the same number of holes are not distinguishuable.

\begin{figure}[ht]
\centering\includegraphics[height=3cm,width=8.cm]{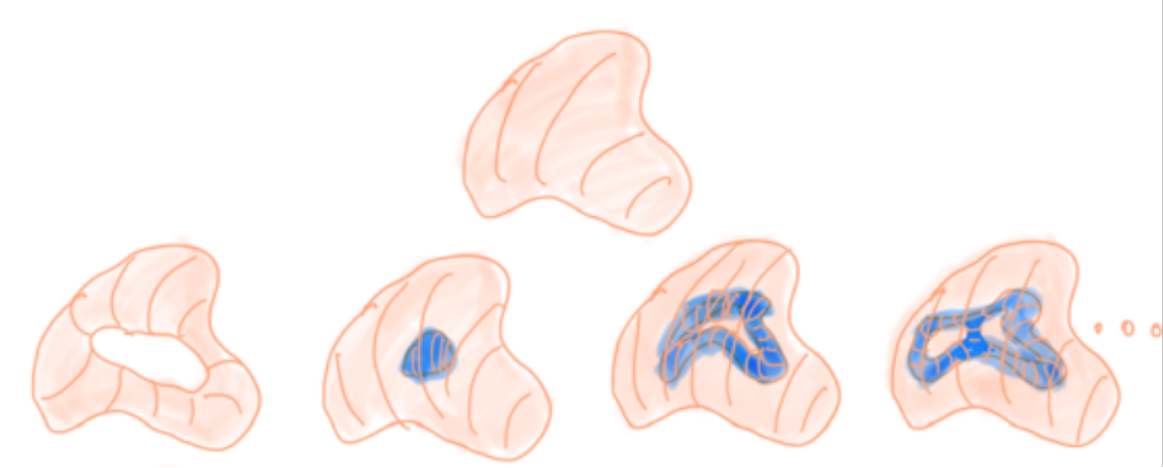}
\caption{Examples of basic topological spaces in three dimensions. {\em Top image}: a simply connected solid space. The brown colour shows solid regions enclosed by a positively oriented surface. {\em Bottom row, left to right}: A connected region having a one-dimensional hole or handle shown by the white region. The second one is a connected region having one two-dimensional hole or cavity within shown in blue.  The third one shows a connected region with a doubly connected cavity or void within, shown in blue. The fourth one  shows a connected region with a triply connected cavity or void within, shown in blue. }
\label{fig:3dspaces}
\end{figure}

{\em Betti numbers for three dimensional topological spaces}: Let us consider three-dimensional topological spaces with boundary. In this case, holes can be of two types - either one or two dimensional. A one-dimensional hole represents the equivalence class of loops on a boundary surface that cannot be shrunk to a point.  A two-dimensional hole is a cavity or void is a 2-dimensional hole, fully enclosed within a 2-dimensional surface or shell.

Fig.~\ref{fig:3dspaces} shows examples of some basic 3D topological spaces having distinct topological features. The top image shows a simply connected solid region. The bottom row, from left to right, shows a connected region with a handle, a connected region with a spherical cavity, a connected region with a tubular cavity, and a connected region with a triply connected cavity. The cavities are shown in blue. 
We may denote these five topological spaces as $Q_{ij_0j_1j_2}$, where $i$ indexes the number of handles, $j_0$ the number of spherical cavities, $j_1$ the number of tubular or doubly connected cavities, and $j_2$ the number of triply connected cavities. For our purpose here topological spaces indexed by the same $ij_0j_1j_2$ are  not distinguished (though they may be further distinguished by other topological properties  such as links and knots, which we ignore).

 The Euler characteristic is given by the alternating sum of the Betti numbers, as,
 \be\chi=b_0-b_1+b_2. \ee
$\chi$ is related to the genus, $g$, by the formula $\chi=2-2g$. For a smooth closed surface, the Gauss-Bonnet theorem relates $g$ to the integral of the Gaussian curvature, $K$, over the surface, as, $\int K \dA = 4\pi(1-g)$. 
The sum of Betti numbers is
\be
b_{\rm sum}= b_0+b_1+b_2. 
\ee
The values of $b_0$, $b_1$, $b_2$, $\chi$ and $\bs$ for each $Q_{ij_0j_1j_2}$ in Fig.~\ref{fig:3dspaces} are given in the following table.\\
\vskip .1cm
\begin{tabular}{p{2cm}p{1cm}p{1cm}p{1cm}p{1cm}p{1cm}}
  \hline
  Space    & \ $b_0$ & \ $b_1$  & \ $b_2$  & \ $\chi$ & \ $\bs$ \\
  \hline
  $Q_{0000}$  &  \ 1 & \ 0 &  \ 0  & \ 1 & \ \ 1\\
  $Q_{1000}$  &  \ 1 & \  1 & \   0  & \  0 & \ \  2\\
  $Q_{0100}$  &  \  1 & \ 0 &  \  1  & \  2 &  \ \ 2\\
  $Q_{0010}$  &  \  1 &  \ 1 &  \  1  & \  1 &  \ \ 3\\
  $Q_{0001}$  &  \  1 &  \ 2 &  \  1  & \  0 & \ \ 4\\
  \hline
  \end{tabular}
  \vskip .3cm

In the next two sections we will be focusing on Betti numbers, Euler characteristic and sum of Betti numbers of random fields in two and three dimensions. For the sake of brevity, we will sometimes use the collective term {\em topological statistics} when we refer to them. 

 \section{Betti numbers, Euler characteristic and sum of Betti numbers for random fields in two dimensions}
\label{sec:s4}

We now focus on random fields on finite two dimensional (2D) smooth manifolds. We first introduce a new topological basis, and express Betti numbers, Euler characteristic and sum of Betti numbers of excursion sets in terms of the new basis. Then we examine their statistical properties that follow from the new representation. We consider fields having zero mean and support in the range $-\infty$ to $\infty$ so as to simplify the discussion.

\subsection{A new topological basis and representations of topological statistics}
\label{sec:s4a}

Let $Q_j$ denote a $j$-multiply connected subset of $M$. $Q_j$ is a connected component that has $j$ number of holes. 
Examples of $Q_j$ upto $j=2$ are shown in Figure~\ref{fig:2dspaces}.

\begin{figure}[b]
  \centering {\large{$\mathcal Q^{(\nu=-1)}$ \hskip 1.6cm $\mathcal Q^{(\nu=0)}$ \hskip 1.6cm $\mathcal Q^{(\nu=1)}$}}\\
    \vskip .1cm
    \centering\includegraphics[height=2.5cm,width=2.5cm]{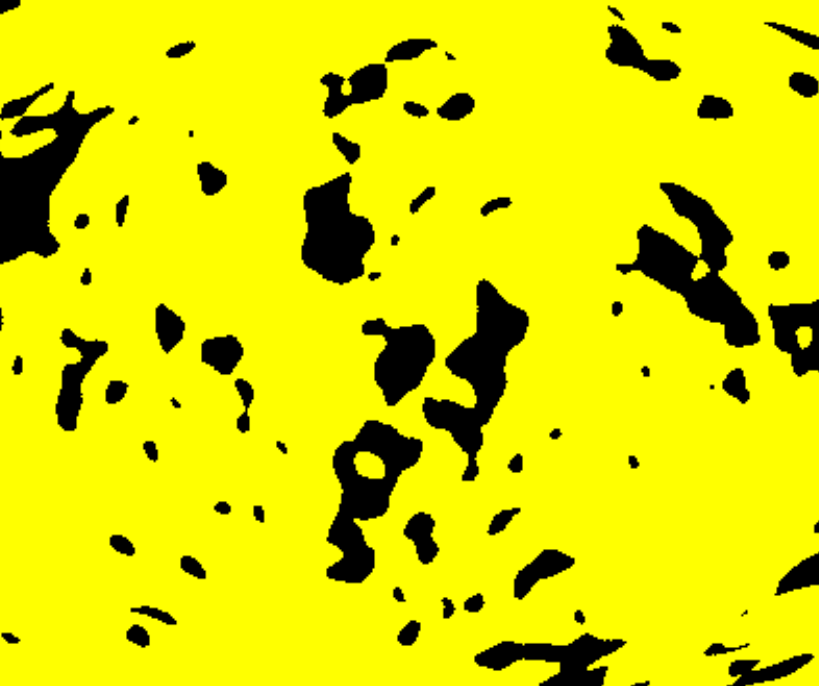} \hskip .3cm
    \centering\includegraphics[height=2.5cm,width=2.5cm]{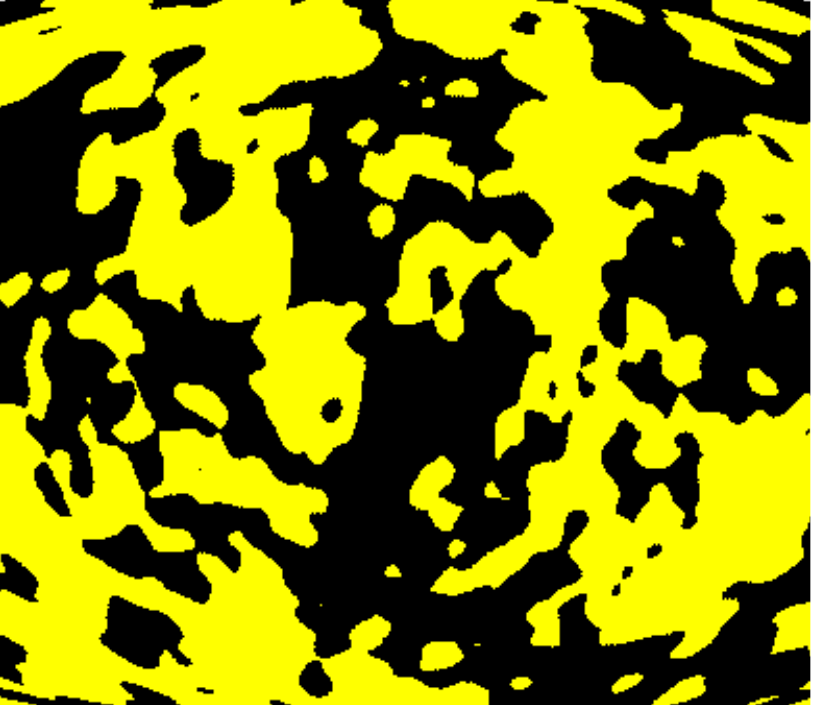} \hskip .3cm
    \centering\includegraphics[height=2.5cm,width=2.5cm]{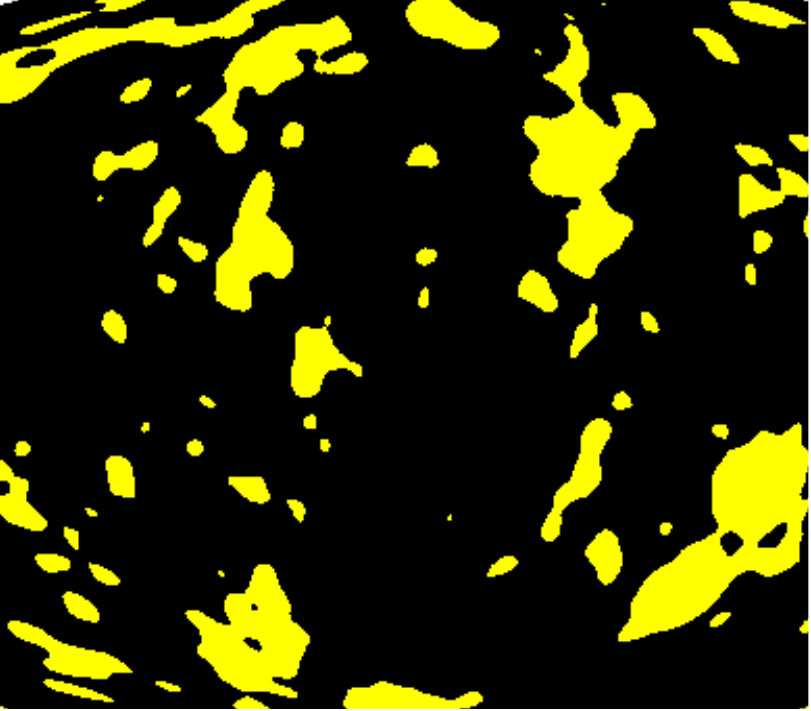}
    \caption{Examples of excursion sets (yellow regions) at three threshold levels for a normalized Gaussian field.}
\label{fig:Qj_2d}
\end{figure}

$Q_j$'s constitute building blocks of excursion sets. To make this statement precise, let $B$ denote the set $\{ Q_0, Q_1,Q_2,...,Q_j,...,Q_{j_{\rm max}} \}$. The maximum index $j_{\rm max}$  is set by the parameter $q$ given by Eq.~(\ref{eq:q}) since it roughly gives the maximum number of structures that can fit in $M$. It is conceivable that there can be configurations of holes whose shapes have high isoperimetric ratios, in which case $j_{\rm max}$ can be larger than  $q$. However, as long as the isoperimetric ratio  is finite $j_{\rm max}$ is finite, and its precise value does not matter for our arguments. 
$B$ forms a topological basis in the sense that an excursion set  $\mathcal{Q}^{(\nu)}$ can be expressed as a {\em linear} union of elements of $B$, as,
\bea
\hspace{-.2cm} \mathcal{Q}^{(\nu)} &=& \UU_{j=0}^{j_{\rm max}} \ m_{j}Q_{j} \nn\\
&=& m_0Q_0 \UU m_1Q_1 \UU m_2Q_2\UU\ldots  \UU m_{j_{\rm max}}Q_{j_{\rm max}}.
\label{eq:Qunion}
\eea
The coefficients $m_j$ are non-negative integers representing the number of copies of $Q_j$ included in $\mathcal{Q}^{(\nu)}$.  $m_j$'s vary with $\nu$ (not explicitly written). Visual demonstration of  Eq.~(\ref{eq:Qunion}) can be seen in Fig.~\ref{fig:Qj_2d} where we show three examples of  $\mathcal{Q}^{(\nu)}$  for a normalized Gaussian field, at threshold values  $\nu =-1,0,1$. Yellow represents connected regions while black represents holes.

Eq.~(\ref{eq:Qunion})  tells us that the randomness of excursion sets are now encoded in $m_j$'s.  Hence, the
crux of understanding the random topology of  $\mathcal{Q}^{(\nu)}$  lies in being able to calculate/guess/estimate the $m_j$'s. 
 Further, the statistical properties of the field should follow from their properties, and vice versa. A schematic flow of statistical properties is shown below.
   \be
    f(\vx) \longleftrightarrow \underset{-\infty< \nu< \infty}{\mathcal{Q}^{(\nu)}}   \longleftrightarrow \ \bigg\{ m_0(\nu),\  m_1(\nu), \ m_2(\nu),...\bigg\}. \nn 
      \ee

For excursion set $\mathcal{Q}^{(\nu)}$, we can express the topological statistics formally in terms of $m_j$'s as
\bea
b_0(\nu) &=& \sum_{j=0}^{j_{\rm max}(\nu)} m_j(\nu), \label{eq:bc}\\
b_1(\nu) &=& \sum_{j=1}^{j_{\rm max}(\nu)}j m_j(\nu),  \label{eq:bv}\\
\chi(\nu) &=& \sum_{j=0}^{j_{\rm max}(\nu)} (1-j)m_j(\nu) \nn\\   
&=& m_0(\nu) - \sum_{j=2}^{j_{\rm max}(\nu)} (j-1)m_j(\nu),  \label{eq:chi}\\
b_{\rm sum}(\nu) &=& \sum_{j=0}^{j_{\rm max}(\nu)}(1+j) m_j(\nu).  \label{eq:bs}
\eea
In the above $j_{\rm max}$ depends on $\nu$ because the maximum $j$ for which $m_{j}$ will be nonzero will depend on $\nu$.

$\chi$ does not contain $m_1$ term since its contribution is zero. Further, for $\chi=0$,  we have $m_0= \sum_{j=2}^{j_{\rm max}} j m_j(\nu)$. For $\chi>0$,  we have $m_0> \sum_{j=2}^{j_{\rm max}} j m_j(\nu) $, and so simply connected components outnumber the sum of all multiply connected components. For $\chi<0$  we have the opposite situation where the sum of all multiply connected components outnumber simply connected components.

Eqs.~(\ref{eq:bc}) and (\ref{eq:bv}) tell us that $b_0,$ and $b_1$  are sums of discrete random variables. So in general, their PDFs will be determined by the PDFs of $m_j$'s, which in turn are  determined by the statistical nature of the field, and vice versa. The properties of  $\chi$ and $\bs$ then naturally follow as the alternating sum and sum of  $b_0,$ and $b_1$ .

 It is interesting to  note that 
$b_0,\ b_1$, $\chi$ and $b_{\rm sum}$ can be expressed in terms of a generating function $h(\a)$, defined as,  
\be
h(\a)= \sum_{j=0}^{j_{\rm max}} m_j\,e^{-j\a},
\ee
where $\a$ is some real variable. This expression is like a discrete form of the Laplace transform  with $j$ being conjugate to $\a$. The difference, however, is that $m_j$ are integers rather than real variables.  
Then we can write
\bea
b_0(\nu) &=&  \lim\limits_{\a\to 0} h(\a)\\
b_1(\nu) &=& - \lim\limits_{\a\to 0}  \frac{dh}{d\a}.
\eea
and
\bea
\chi(\nu) &=&   h(\a)|_{\a=0}  +  \frac{dh}{d\a}\bigg|_{\a=0},\\
b_{\rm sum}(\nu) &=&   h(\a)|_{\a=0}  -  \frac{dh}{d\a}\bigg|_{\a=0}.
\eea

\subsection{Statistical nature of the topological statistics}
\label{sec:s4b}

To fully exploit the new basis for deeper understanding of the statistical nature of the topological statistics we require information about $m_j$'s. While a full theoretical exploration of the properties of $m_j$'s is beyond the scope of this paper, we proceed by making the simplifying assumption that at each threshold $m_j$'s are statistically independent. 
Under this condition their joint probability distribution function gets factorised as
\be
P\left(m_0(\nu),m_1(\nu),...\right)= P(m_0(\nu))P(m_1(\nu))....
\label{eq:mj_jpdf}
\ee
Then the ensemble expectations of $b_0,b_1, \chi$ and $\bs$ are 
\bea
\big\la b_0(\nu)\big\ra &=& \sum_{j=0}^{j_{\rm max}} \big\la m_j(\nu)\big\ra,\\
\big\la b_1(\nu)\big\ra &=& \sum_{j=1}^{j_{\rm max}}j \big\la m_j(\nu)\big\ra,\\
\big\la \chi(\nu)\big\ra &=& \sum_{j=0}^{j_{\rm max}} (1-j)\big\la m_j(\nu)\big\ra,\\
\big\la b_{\rm sum}(\nu)\big\ra &=& \sum_{j=0}^{j_{\rm max}} (1+j)\big\la m_j(\nu)\big\ra.  
\label{eq:gmj}
\eea
Note that $\nu$ dependence of $j_{\rm max}$ is understood. 
Their variances are  
\bea
\s_{b_0}^2(\nu) &=& \sum_{j=0}^{j_{\rm max}} \big\la m_j^2(\nu)\big\ra - \sum_{j,k=0}^{j_{\rm max}} \big\la m_j(\nu)\big\ra \big\la m_k(\nu)\big\ra ,\label{eq:sb0}\\
\s_{b_1}^2(\nu) &=& \sum_{j=1}^{j_{\rm max}}j^2 \big\la m_j^2(\nu)\big\ra -  \sum_{j,k=0}^{j_{\rm max}} jk\big\la m_j(\nu)\big\ra \big\la m_k(\nu)\big\ra , \nn \label{eq:sb1} \\ \\
\s_{\chi}^2(\nu) &=& \sum_{j=0}^{j_{\rm max}} (1-j)^2\big\la m_j^2(\nu)\big\ra \nn\\
&& -  \sum_{j,k=0}^{j_{\rm max}} (1-j)(1-k)\big\la m_j(\nu)\big\ra \big\la m_k(\nu)\big\ra, \label{eq:schi}\\
\s_{b_{\rm sum}}^2(\nu) &=& \sum_{j=0}^{j_{\rm max}} (1+j)^2\big\la m_j^2(\nu)\big\ra \nn\\
&& -  \sum_{j,k=0}^{j_{\rm max}} (1+j)(1+k)\big\la m_j(\nu)\big\ra \big\la m_k(\nu)\big\ra . \label{eq:sbsum}
\eea

Eqs.~(\ref{eq:schi}) and (\ref{eq:sbsum}) can also be written as
\bea
  \sigma^2_{\chi} &=& \sigma^2_{b_0} + \sigma^2_{b_1} - 2{\rm Cov}\left(b_0,b_1\right), \label{eq:schi2}\\
  \sigma^2_{b_{\rm sum}} &=& \sigma^2_{b_0} + \sigma^2_{b_1} + 2{\rm Cov}\left(b_0,b_1\right),  \label{eq:sbsum2}
\eea
where the last term on the right hand sides is the covariance of $b_0$ and $b_1$ given by 
\bea
    {\rm Cov}\left(b_0,b_1\right) &=& \langle b_0b_1 \rangle -  \langle b_0 \rangle  \langle b_1 \rangle\nn\\
    &=& \sum_j j\langle m_j^2\rangle  - \sum_{jk} j\langle m_j\rangle  \langle m_k\rangle. \label{eq:Covb0b1}     
\eea
Intuitively, we may expect $b_0$ and $b_1$ to be anti-correlated, with a negative  $ {\rm Cov}\left(b_0,b_1\right)$. This expectation arises from the observation that when the number of connected components is large, the number of holes tends to be small, and vice versa. This behavior is confirmed by our numerical calculations for Gaussian fields presented in the next subsection.

\subsubsection{Comparison of standard deviations of the topological statistics for Gaussian fields}
\label{sec:s4b1}

\begin{figure*}[t]
    Ensemble expectations \hskip 2.4cm Standard deviations \hskip 3.9cm Cov($b_0,b_1$)\\
  \centering\includegraphics[width=.31\linewidth]{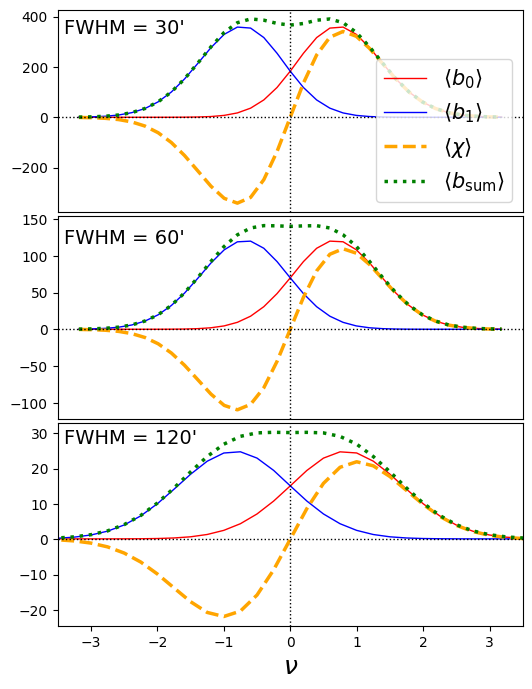} \hskip .3cm
  \centering\includegraphics[width=.3\linewidth]{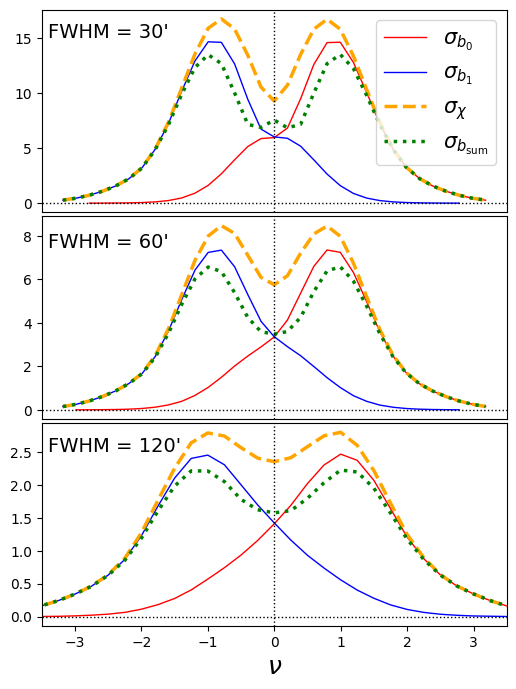} \hskip .3cm
   \centering\includegraphics[width=.31\linewidth]{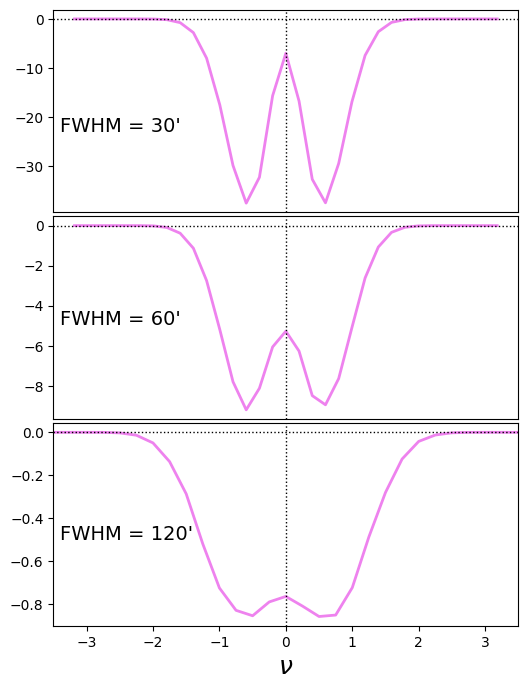}
\caption{Ensemble expectations (left) and standard deviations (middle) of $b_0,\,b_1,\,\chi$ and $b_{\rm sum}$ computed numerically from $10^4$ CMB temperature maps smoothed with Gaussian kernels having FWHM $30'$, $60'$ and $120'$. The panels in the right column show the covariance between $b_0$ and $b_1$.}
\label{fig:mean_sigma}
\end{figure*}
Before using the topological basis further we examine numerically computed ensemble expectations and standard deviations of $b_0$, $b_1$, $\chi$ and $b_{\rm sum}$.
  For this purpose, we simulate  $10^4$  Gaussian CMB temperature fluctuation maps whose pixel resolution is set by the Healpix parameter $N_{\rm side}=512$. The maps are then smoothed with Gaussian kernel with full width at half maximum (FWHM)  $30'$, $60'$ and $120'$. The four statistics are then computed for the smoothed maps.

The left  column of Fig.~\ref{fig:mean_sigma} show the ensemble expectations of $b_0$, $b_1$, $\chi$ and $b_{\rm sum}$, {\em per unit area} (note that we use the same symbols), for three smoothing scales. For Gaussian fields,  $\langle \chi\rangle$  per unit area has a closed form analytic formula~\cite{Adler:1981,Tomita:1986} given by
\be \langle \chi(\nu)\rangle =A_2 \nu e^{-\nu^2/2},
\label{eq:chig}
\ee
where  $A_2=1/(4\sqrt{2}\pi^{3/2} r_c^2)$ is the amplitude. 
Analytic formula for Betti numbers are not known. But they are well studied numerically~\cite{Park:2013,Pranav:2017}. The plots tell us that in negative thresholds the topology of excursion sets is dominated by holes, while for positive thresholds connected regions dominate.  In the regime $\nu\to \infty$, we have $\chi,\bs\to b_0$. In contrast, for $\nu\to -\infty$, we have $\chi\to -b_1$ and $\bs\to b_1$. The roles of connected regions and holes get reversed under a field transformation $f\to -f$.

The middle  column of Fig.~\ref{fig:mean_sigma} shows the standard deviations of the four statistics for the three smoothing scales. What is interesting to note is that we get
  \bea
  \s_{\chi} &>& \sqrt{\sigma_{b_0}^2 + \s_{b_1}^2},\\
  \s_{\bs} &<& \sqrt{\sigma_{b_0}^2 + \s_{b_1}^2}.
  \eea
This is because the variances of $\chi$ and $b_{\rm sum}$ relative to the variances of $b_0$ and $b_1$ are determined by ${\rm Cov}\left(b_0,b_1\right)$ (see Eqs.~(\ref{eq:schi2}) and (\ref{eq:sbsum2})) which we obtain to be negative.   
The right column of Fig.~\ref{fig:mean_sigma} shows  ${\rm Cov}\left(b_0,b_1\right)$ for the three smoothing scales. Therefore, $b_0$ and $b_1$ are anti-correlated, as anticipated earlier. The shape of the $\nu$ dependence varies with FWHM. 

The negative sign of ${\rm Cov}\left(b_0,b_1\right)$ further implies the following,
\be
 \sum_j j \bigg( \langle m_j^2\rangle  - \langle m_j\rangle \sum_{k} \langle m_k\rangle \bigg) < 0, 
 \ee
 which can be further expressed as
 \be
 \sum_j j  \langle m_j\rangle  \bigg(\langle m_j\rangle  + \frac{\s^2_{m_j}}{\langle m_j\rangle}  - \sum_{k} \langle m_k\rangle \bigg) < 0.  \label{eq:con}
 \ee
 There are two possibilities for satisfying this inequality. The first is that the following strict condition is satisfied for every $j$,
  \be
\langle m_j\rangle  + \frac{\s^2_{m_j}}{\langle m_j\rangle}<  \sum_{k} \langle m_k\rangle. \label{eq:strictcon}
\ee
Since in general $\langle m_j\rangle <  \sum_{k} \langle m_k\rangle $, the above condition requires the variance of $m_j$ to be small relative to its mean. 
The second possibility is that the inequality Eq.~\ref{eq:strictcon} may be violated for some $j$'s, provided the total sum over over $j$ remains negative.

The discussion above uses numerical computation to shed light on the properties of $m_j$'s for Gaussian fields. It is worth mentioning that 
the relatively lower standard deviations of $b_{\rm sum}$ compared to the other statistics can potentially result in tighter constraints when carrying out cosmological parameter inference.   

\subsubsection{Central limit theorem for topological statistics} 
\label{sec:s4b2}

We are interested in the behaviour of $b_0, b_1$, $\chi$ and $\bs$ in the limit of large $j_{\rm max}$. For this purpose let us recapitulate the central limit theorem (CLT) for statistically independent discrete random variables~\cite{Lindeberg:1922,Feller:1945}.  Let $X_{i}$ be $n$ independent discrete random variables that are  not necessarily identically distributed. Let each $X_i$ have mean $\mu_i$ and variance $\s_i^2$. Let $X = \sum _{i=1}^{n} X_i$, and $s_{n}^{2}= \sum _{i=1}^n\sigma_i^2$. 
Assume that $s_n \to \infty$ as ${n\to \infty}$. This is satisfied unless $\s_i \to 0$ as $i\to\infty$. Then,  $X$ is asymptotically normally distributed if there exists a constant $c$ such that: 
\be
|X_i| \le c, \ {\rm for \ all} \ i.
\label{eq:LF}
\ee
For sequences of random variables that are block-sum independent (see e.g.~\cite{Dvoretzky:1970,Fleermann:2022}), we  can write $X$,
\bea X = X_1+ X_2+ \ldots +X_n = Y_1 + Y_2+ .. +Y_M,
\eea
where $Y_i$'s are block-sums of $X_i$'s that are dependent in a block of size $k$, such that $M k =n$, and $k $ is greater tan the correlation length. Then treating $Y_i$'s as independent random variables, one can appeal to the CLT.

Now we examine if the CLT holds for $b_0, b_1$, $\chi$ and $\bs$  without knowing the explicit form of PDFs of $m_j$'s. For this we require $j_{\rm max}\to \infty$, or $q\to \infty$,  which further depends on the form of the power spectrum as summarised in table \ref{tab:t1}. 
Assuming that the power spectrum is such that $q\to\infty$, let us consider the $m_j$'s to be at least block-sum statistically independent. 
Then ${\widetilde m}_j\equiv jm_j$'s are also statistically independent. 

{\em CLT for $b_0$}: Let $X_i = m_{i-1}$  and $n=j_{\rm max}-1$.  Then $X=b_0$. 
We need to check whether the condition set by Eq. (\ref{eq:LF}) is satisfied by $X_i$ or not. For a given resolution set by $q$, the variable having the largest (finite) support is $m_0$, and it satisfies the bound $m_0< j_{\rm max}$.  Hence Eq. (\ref{eq:LF}) is satisfied. Therefore, $b_0$ must be asymptotically  Gaussian for ${j_{\rm max}\to\infty}$.

{\em CLT for $b_1$}: Next, let $X_i = \widetilde m_{i-1}$  and $n=j_{\rm max}-1$.  Then $X=b_1$. 
In this case, for resolution  $q$, the variable having the largest (finite) support is $m_{j_{\rm max}}$, and it satisfies the bound $m_{j_{\rm max}}\le j_{\rm max}$. Therefore, $b_1$ must  also be asymptotically  Gaussian for ${j_{\rm max}\to\infty}$.

{\em CLT for $\chi$ and $\bs$}: For $\chi$, let $m'_j =(1-j) m_j$. Then $\chi$ can be written as
\be
\chi = \sum_{j=0}^{j_{\rm max}} m'_j.
\ee
In this case also $|m'_j|$ are bounded, with the largest element being determined by whether $b_0>b_1$ or $b_0<b_1$. Therefore, $\chi$ is also asymptotically  Gaussian if ${j_{\rm max}\to\infty}$.
Finally, for $\bs$,  let $m''_j =(1+j) m_j$. Then $\bs$ can be written as
\be
\bs = \sum_{j=0}^{j_{\rm max}} m''_j.
\ee
with $|m''_j|$  bounded by the largest element determined by whether $b_0>b_1$ or $b_0<b_1$. Therefore, $\bs$ is also asymptotically  Gaussian if ${j_{\rm max}\to\infty}$.

To summarize, we have shown that $b_0$, $b_1$, $\chi$ and $\bs$  at each $\nu$ tend to Gaussian random variables, provided $q\to\infty$. {\em We stress  that due to the dependence of $j_{\rm max}$  on $\nu$ the convergence of the topological statistics to Gaussian nature depends on $\nu$.}  

\subsubsection{Modeling $m_j$'s as Binomial distributed variables}
\label{sec:s4b3}

In practical applications the limit $q \to \infty$ discussed in the previous subsection is not of much use since observed data covers finite spatial extent and finite resolution. 
We now take a different approach and focus on the nature of $m_j$'s. 

Consider a box containing balls of distinct colours. Let each element of the topological basis, $Q_j$, be represented by a distinct colour ball.  
{\em An excursion set at $\nu$ can be considered as a realization of $m_j$ number of each distinct ball in $N$ trials of drawing balls from the box.} {Note that the size  and shape variations of the balls are ignored in this statement since we are focusing only  on their topological properties.} Then, each $m_j$ is a Binomial variable which can take value between 0 to $N$.  Let $p_j$ denote the success probability that the ball is type $Q_j$. 
In $N$ trials the probability that a $Q_j$ type ball will be picked $k$ times is given by the Binomial distribution, 
\be
P(m_j= k;p_j,N)={}^NC_k\,p_j^k\,(1-p_j)^{N-k}.
\ee
How $p_j$ varies with  $j$, and how both $p_j$ and $N$  depend on $\nu$ will be dictated by the nature of the field.   
The mean and variance of $m_j$ are then given by
\bea
\langle m_j \rangle &=& p_jN, \\
\label{eq:mmj}
\sigma_{m_j}^2 &=& Np_j(1-p_j).
\label{eq:smj}
\eea

Since $b_0$ is given by the sum of $m_j$'s, its PDF is given by the product of their  PDFs, 
\be
P(b_0) = P(m_0)P(m_1)\dots,
\label{eq:b0pdf}
\ee
provided $m_j$'s are independent of each other.  $P(b_0)$  is, in general, not a Binomial distribution since $p_j$'s need not be the same for different $j$. The mean and variance of $b_0$  are, 
\bea
\langle b_0(\nu) \rangle &=& \sum_{j=0}^{j_{\rm max}}p_j(\nu)N(\nu), \\
\label{eq:mb0}
\sigma_{b_0}^2(\nu) &=& \sum_{j=0}^{j_{\rm max}}p_j(\nu)(1-p_j(\nu))N(\nu).
\label{eq:sb0_bino}
\eea

For $b_1$, if $m_j$ is binomial distributed, then $\widetilde m_j\equiv jm_j$ is also binomial distributed, with probability $\widetilde p_j$:  
\be
P(\widetilde m_j= k;\widetilde p_j,N)={}^NC_j\,\widetilde p_j^{k}\,(1-\widetilde p_j)^{N-k}.
\ee
Therefore,
\be
P(b_1) = P(\widetilde m_1)P(\widetilde m_2)\dots.
\label{eq:b1pdf}
\ee
The mean and variance of $b_1$  are, 
\bea
\langle b_1(\nu) \rangle &=& \sum_{j=1}^{j_{\rm max}}\widetilde p_j(\nu)N(\nu), \\
\label{eq:mb1}
\sigma_{b_1}^2(\nu) &=& \sum_{j=1}^{j_{\rm max}}\widetilde p_j(\nu)(1-\widetilde p_j(\nu))N(\nu).
\label{eq:sb1_bino}
\eea

The PDFs of $\chi$ and $\bs$ must then follow from the PDFs of $b_0$ and $b_1$. 
We stress that while the Binomial nature of the $m_j$'s as modelled above is general and does not depend on the nature of the field,
the specific  forms of their PDFs depend on $p_j(\nu)$'s and $N(\nu)$, which must be determined by the nature of the field. 

\subsubsection{Testing the Binomial modeling for the case of Gaussian fields}
\label{sec:s4b4}

Since  $p_j(\nu)$'s and $N(\nu)$ are not known in general, we now adopt a semi-numerical approach.     
Specifically we focus on  Gaussian random fields, for which the closed form analytic formula for $\langle\chi\rangle$ is available. This allows us to test the validity of the Binomial model in relevant threshold regimes. 

As mentioned earlier, a consequence of  Eq.~(\ref{eq:chi}) is that  $m_0$ outnumber the sum over all other $m_j$'s for $\chi>0$. For Gaussian fields with zero mean this holds for $\nu> 0$.  As $\nu$ increases in the positive direction we can approximately write
\bea
b_0(\nu) &\sim&  m_0(\nu),\\
b_1 &\sim& 0,\\
\chi(\nu) &\sim &  m_0(\nu),\\
b_{\rm sum}(\nu) &\sim&  m_0(\nu).
\eea
The accuracy of the approximations increases as $\nu$ increases. Using the right hand side of  Eq.~(\ref{eq:chig})
and equating it  to $ \langle b_0(\nu) \rangle \sim p_0(\nu)N(\nu)$, we get 
\be
p_0(\nu) = \frac{A_2}{N(\nu)} \nu e^{-\nu^2/2}.
\label{eq:p+}
\ee
Inserting the above in Eq.~(\ref{eq:smj}) for $j=0$, we get,
\bea
\sigma_{m_0}^2(\nu) &\sim&\sigma_\chi^2\sim  A_2 \nu e^{-\nu^2/2}\left(1-\frac{A_2}{N(\nu)} \nu e^{-\nu^2/2} \right). \quad
\label{eq:s+}
\eea
Then, using numerically computed $\sigma_{\chi}$ and equating with the right hand side of Eq.~(\ref{eq:s+}), we can obtain $N(\nu)$.

\vskip .2cm  
For $\nu < 0$, we have $\chi<0$, and the sum of all $m_j$'s with $j\ge 1$ outnumber $m_0$.  As $\nu$ decreases  we can approximate as  $\chi \to -b_1$, $b_{\rm sum} \to b_1$. In this regime the excursion set tends to consist of a single multiply connected region. So we have, 
\bea
b_0(\nu) &\sim&  m_{j_{\rm max}} \sim 1,\\
b_1(\nu) &\sim& j_{\rm max}(\nu) m_{j_{\rm jmax}}(\nu) \sim j_{\rm max},\\
\chi(\nu) &\sim &  - j_{\rm max}(\nu) m_{j_{\rm max}}(\nu) \sim - j_{\rm max},\\  
b_{\rm sum}(\nu) &\sim& j_{\rm max}(\nu) m_{j_{\rm jmax}}(\nu) \sim j_{\rm max},
\eea
Then, using Eq.~(\ref{eq:chig}) we get
\be
p_{j_{\rm max}}(\nu) \sim -\frac{A_2}{N} \nu e^{-\nu^2/2},
\label{eq:p-}
\ee
and
\bea
\sigma_{m_{j_{\rm max}}}^2\sim \s_{\chi}^2 &\sim & -A_2 \nu e^{-\nu^2/2}\left(1+\frac{A_2}{N(\nu)} \nu e^{-\nu^2/2} \right). \quad\ \ 
\label{eq:s-}
\eea
Again by using numerically computed $\sigma_{\chi}$ and equating with the right hand side of Eq.~(\ref{eq:s-}), we can obtain $N(\nu)$.

For Gaussian random fields, Betti numbers of high and low density thresholds are related as $b_0(\nu) = b_1(-\nu)$. This is a reflection of Alexander duality~\cite{alexander:1923,hatcher:2002} which relates the topology of high-threshold excursion regions to the topology of low-threshold regions. 

We now test how well numerically computed PDFs of $b_0$, $b_1$, $\chi$ and $b_{\rm sum}$ agree with the Binomial modeling of $m_j$'s, using the same CMB temperature simulations as in  Fig.~\ref{fig:mean_sigma}. To do this, we  calculate semi-analytic approximate Binomial PDFs for each statistic using numerically computed $\s_{\chi}$ and analytic formula for $\langle\chi\rangle$. 
We then compute numerical PDFs at each threshold and compare these results with the Binomial predictions. We expect good agreement between the two, particularly towards  $|\nu|>1$. 
Additionally, we compare these results with Gaussian PDFs constructed from the ensemble means and standard deviations of each statistic at each threshold.

\begin{figure*}[t]
  \centering {\bf Positive $\nu$}\\
  \vskip .2cm
   \centering\includegraphics[scale=.13]{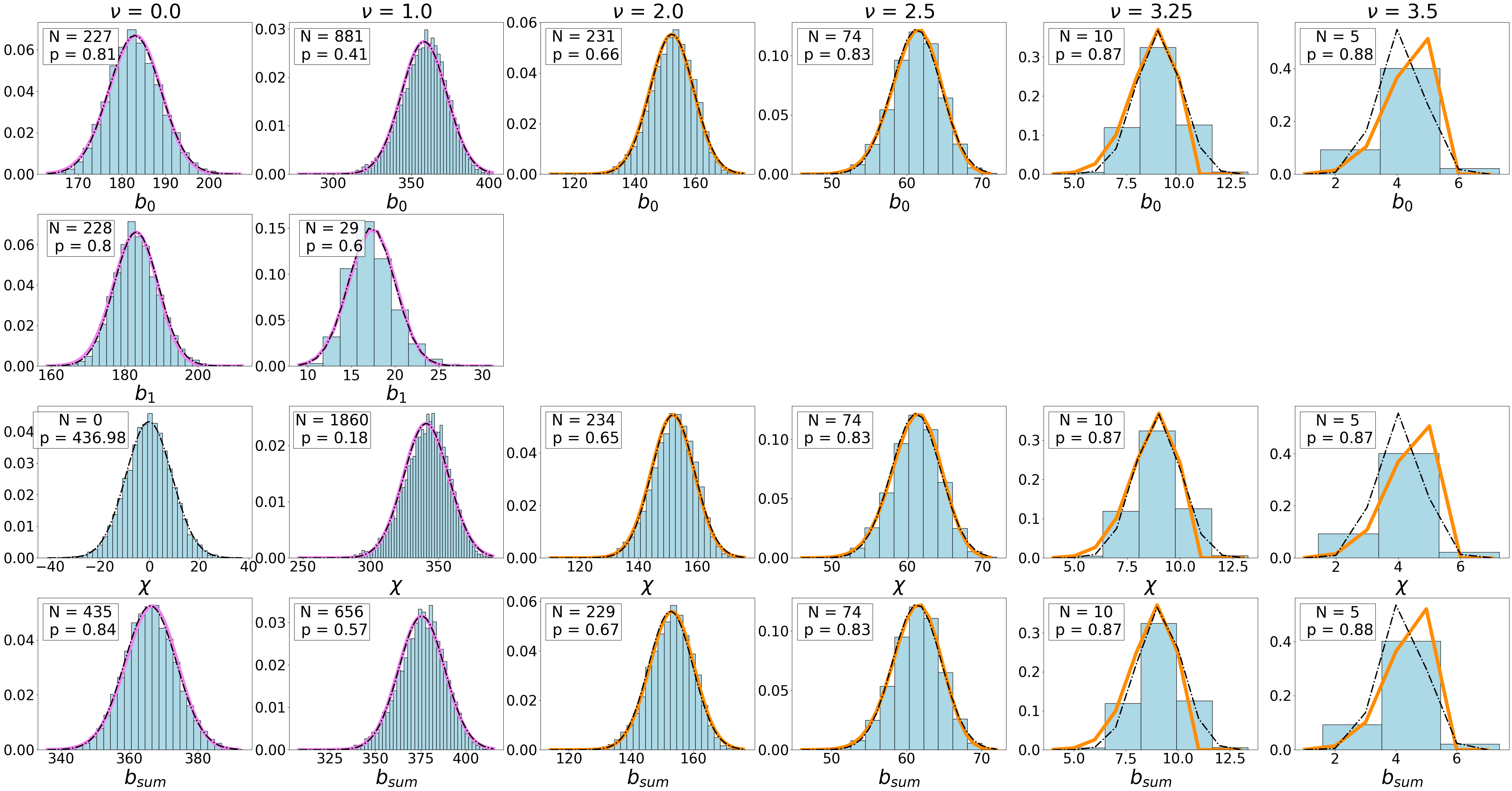}\\
\vskip .3cm
\centering  {\bf Negative $\nu$}\\
\vskip .2cm
\centering\includegraphics[scale=.13]{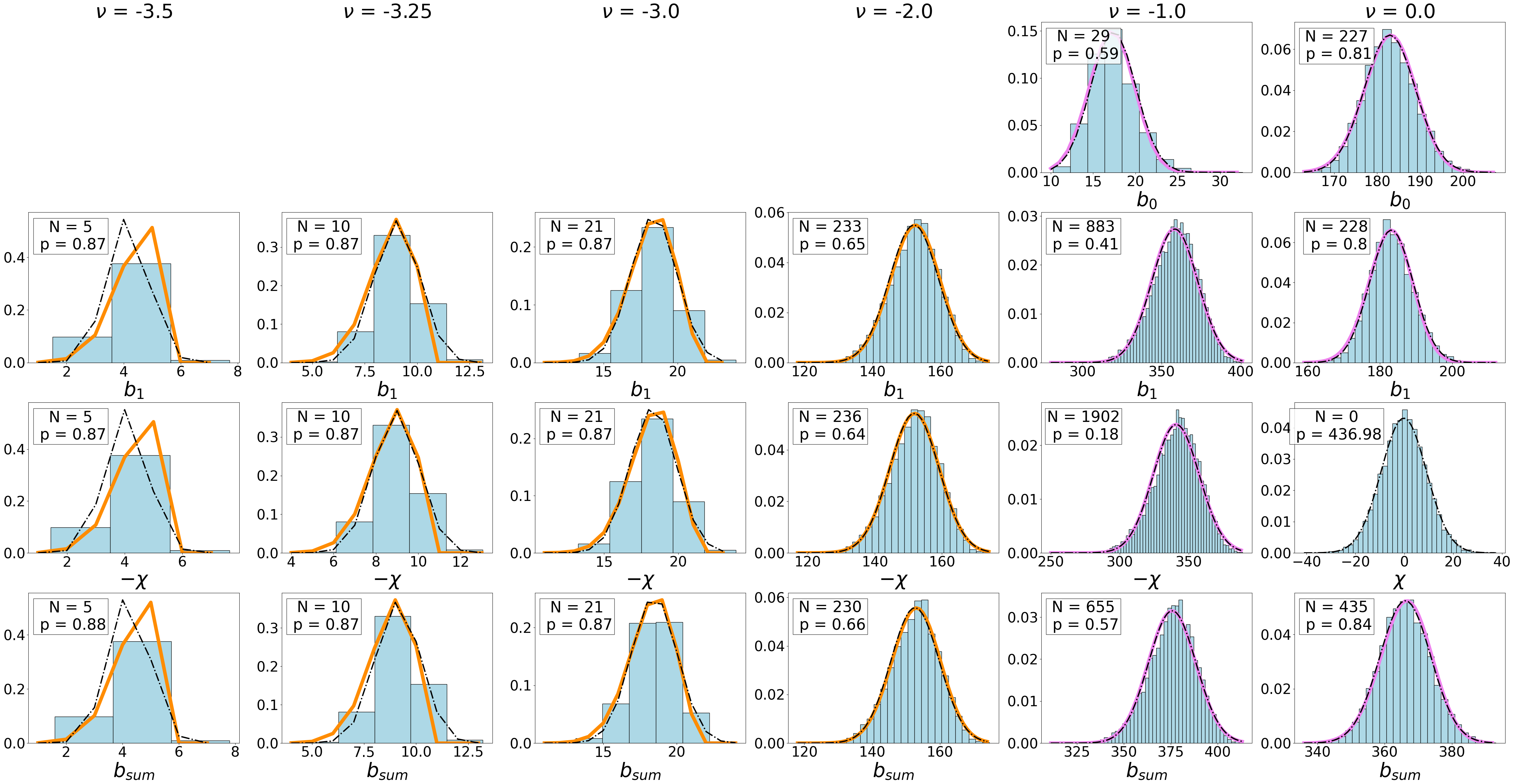}  
\caption{PDFs of $b_0$, $b_1$, $\chi$ and $b_{\rm sum}$ obtained from $10^4$ simulations of Gaussian CMB temperature maps with smoothing scale FWHM $30'$.
  The values of $N,p$ shown are obtained using the ensemble expectation and standard deviation of each statistic, at each $\nu$, as inputs in the Binomial distribution. 
  The solid orange curves are Binomial curves in the regimes of large positive $\nu$ where $b_0\sim\bs\sim\chi\sim m_0$ with $p\sim p_0$, and large negative $\nu$ where $ b_1\sim\bs \sim -\chi \sim j_{\rm max}$ with $p\sim p_{j_{\rm jmax}}$.  The magenta curves in the intermediate regime $-2 < \nu < 2$  are obtained by assuming all four statistics to be Binomial variables and obtaining  $N,p$ using the ensemble expectation and standard deviation of each statistic. 
  The black dashed lines are Gaussian PDFs plotted using the ensemble expectations and standard deviations of $b_0$, $b_1$, $\chi$ and $b_{\rm sum}$.} 
\label{fig:pdfs}
\end{figure*}

The PDFs are shown in Fig.~\ref{fig:pdfs}.  We discuss the results separately for three threshold regimes, as follows.
\begin{itemize}
\item {\em Large positive regime},  $\nu > 1$:  In each panel the values of $N$ and $p$ shown for $b_0$, $\chi$ and $b_{\rm sum}$ are obtained using Eqs.~(\ref{eq:p+}) and (\ref{eq:s+}), with  numerically computed $\s_{\chi}$ as inputs. Note that we dropped the index for $p$, and in this regime $p\sim p_0$. $b_1$ is not shown because it approaches zero in this regime. The solid orange lines represent the corresponding Binomial distributions, which  fit the numerical PDFs well. The black dashed lines show Gaussian PDFs, which generally agree with the Binomial curves except at very large $\nu$. This is expected  because the Binomial distribution approaches a the Gaussian form in the limit of large $N$. For $\nu>3$, $N$ is of order one and $p_0\to 1$, causing the Gaussian approximation to  deviate from the Binomial model, as anticipated. 
\item {\em Large negative regime},  $\nu< -1$:  $N$ and $p$ shown in each panel for $b_1$, $\chi$ and $b_{\rm sum}$ are obtained using Eqs.~(\ref{eq:p-}) and (\ref{eq:s-}), with numerically computed $\s_{\chi}$ as inputs.  In this regime we have $p\sim p_{j_{\rm max}}$. $b_0$ is not shown because it drops to zero in this regime. The solid orange lines are again the corresponding Binomial curves. They fit  the numerical PDFs well. The black dashed lines are Gaussian PDFs. 
They agree well with the Binomial curves for large $N$ and deviates for low $N$, as expected. 
\item {\em Intermediate  regime},   $-2 < \nu < 2$:   In this regime, both $b_0$ and $b_1$ receive contributions from multiple  $m_j$'s, making them sums of Binomial variables, each characterized by  distinct $p_j$ values.  Nevertheless, we {\em assume} $b_0$ and $b_1$ to be single Binomial variables and estimate $N$ and $p$ (shown in the panels) using their numerically computed ensemble means and standard deviations as inputs. The solid magenta lines represent the corresponding Binomial distributions, while the black dashed lines show Gaussian curves. Similarly, Binomial, numerical and Gaussian PDFs are also computed likewise for $\chi$  and $\bs$. At $|\nu|=1$,  the magenta and black lines agree well with the numerical PDFs, although the values of $N$ and $p$ differ between the statistics, as expected.
 
 In the panels for $|\nu|=0$, the numerical, Binomial and Gaussian PDFs for $b_0$ and $b_1$ are all comparable, as expected for Gaussian fields.  For $\chi$, however, we Binomial fit to the numerical PDF is not possible because its expectation value is zero. This leads to unphysical values $N=0$ and $p\gg 1$. The Gaussian curve (shown in black) fits the numerical PDF well. Finally,  the numerical PDF of $b_{\rm sum}$ is well approximated by the Binomial and Gaussian models. Here, the estimated $N$ is roughly twice that for $b_0$ and $b_1$, but the value of $p$ is comparable. 
    
\end{itemize}

{\em Validity of Gaussian approximations}: As discussed in Section~\ref{sec:s1}, assessing whether the PDFs of topological statistics are Gaussian is crucial for cosmological inference. Our results show that the new topological representation, combined with Binomial modeling of coefficients, naturally explains why Gaussian approximations hold across different threshold regimes. 
The approximation is accurate when  $N$ is large, but breaks down when $N$ is order one. 
The amplitudes of the topological statistics as functions of $\nu$ decrease with decrease of resolution. This results in decrease of $N$ with resolution. Therefore, the accuracy of the Gaussian approximation diminishes at lower resolution. 

Note that by modeling $m_j$'s as Binomial variables we are assuming Bernoulli trials without replacement.  It we  model with replacement then $m_j$'s will be Hypergeometric variables. In the limit of large number of trials the Hypergeometrc distribution tends to the Binomial distribution. Hence, the discussion above holds, except at high thresholds.   

\subsection{Counting of states}
\label{sec:s4c}

As a digression we ask how many distinct sequences of $m_j$'s can correspond to a given pair of values $(b_0,b_1)$ by using a combinatorial argument. 
Let us call the sequence $\psi\equiv (m_0,m_1,...,m_{j_{\rm max}})$, for which we get a pair of values $(b_0,b_1)$, as a {\em state}. In general, more than one state can correspond to each $(b_0,b_1)$.  Let the number of distinct states corresponding to a fixed $(b_0,b_1)$ be denoted by $N_{\psi}$. Below we calculate $N_{\psi}$ for three separate cases.
\vskip .2cm
\noindent{{\bf Case 1: $b_0=b_1 =n$.}} In table~\ref{tab:t2} we explicitly write all possible $\psi$ shown by the columns of values of $(m_0,m_1,...,m_n)$  for $b_0=b_1=n$, for $n=2,3,4$. By induction, for arbitrary $n$ we obtain $N_{\psi}=n$.
\renewcommand{\arraystretch}{1.2}
\begin{center}
\begin{table}[ht]
\begin{tabular}{|p{1.8cm}|c|c|c|} \hline
\  & $n= 2$ & $n=3$ & $n=4$  \\ \hline 
{\bf States} & 
\begin{tabular}{cccc}
  $m_0$ & $\rightarrow$  & 1 & 0 \\
  $m_1$ & $\rightarrow$  & 0 & 2 \\
  $m_2$ & $\rightarrow$  & 1 & 0 \\
  &&& \\
  &&& 
  \end{tabular}       & 
\begin{tabular}{ccccc}
  $m_0$ & $\rightarrow$  & 2 & 1 & 0  \\
  $m_1$ & $\rightarrow$  & 0 & 1 & 3 \\
  $m_2$ & $\rightarrow$  & 0 & 1 & 0 \\
  $m_3$ & $\rightarrow$  & 1 & 0 & 0   \\
  &&& &
  \end{tabular}       & 
  \begin{tabular}{cccccc}
  $m_0$ & $\rightarrow$  & 3 & 2 & 1 & 0  \\
  $m_1$ & $\rightarrow$  & 0 & 1 & 2 & 4  \\
  $m_2$ & $\rightarrow$  & 0 & 0 & 1 & 0 \\
  $m_3$ & $\rightarrow$  & 0 & 1 & 0 & 0   \\
  $m_4$ & $\rightarrow$  & 1 & 0 & 0 & 0   
  \end{tabular} \\        
\hline
    $N_{\psi}$ & 2 & 3 & 4 \\ \hline
\end{tabular}
\caption{Number of states for the case $b_0=b_1=n$. For each $n$, each column of values  $m_0, m_1,...$ represents a state.}
\label{tab:t2}
\end{table}
\end{center}

\noindent {\bf Case 2: $b_0> b_1$,  $b_0=n_0, b_1=n_1$.} First, $m_0$ can take values $n_0-k$, where k can be $1,2,...,n_1$. 
This gives $n_1$ possibilities. For each possibility, the problem boils down to: in how many ways can $n_1$ indistinguishable balls (here  holes) be placed in $k$ number of distinguishable boxes (here connected regions), such that each box gets at least one ball. 
This is the so called `stars and bars' problem~(see e.g.~\cite{Feller:1968}). The answer is given by  ${}^{n_1-1}C_{k-1}$. Then summing over all possible values of $m_0$, we get $N_{\psi}$ to be
\be
N_{\psi} = \sum_{k=1}^{n_1} {}^{n_1-1}C_{k-1}.
\label{eq:Ns2}
\ee
Note that $N_{\psi}$ depends only on $n_1$, and $n_0$ does not enter the formula explicitly.

\vskip .2cm
\noindent{{\bf Case 3: $b_0< b_1$,  $b_0=n_0, b_1=n_1$.}} 
 $m_0$ can take values $1,2,...,n_0$. For each value of $m_0=k$, the counting here is:  in how many ways can $n_1$ indistinguishable holes be placed in $k$ number of distinguishable connected regions,  such that each connected region gets at least one hole. The counting is identical to case 2, and  the answer is again ${}^{n_1-1}C_{k-1}$. Again summing over all possible values of $m_0$, $N_{\psi} $ is given by
\be
N_{\psi} = \sum_{k=1}^{n_0} {}^{n_1-1}C_{k-1}.
\label{eq:Ns3}
\ee
Eq.~\ref{eq:Ns3} differs from Eq.~\ref{eq:Ns2} only in the limit of the summation. 

 \section{Betti numbers, Euler characteristic and sum of Betti numbers for random fields in three dimensions}
 \label{sec:s5} 
 
 Excursion sets of random fields on three-dimensional (3D) manifolds exhibit significantly more complex topological behavior than their two-dimensional counterparts. A typical excursion set in 3D may comprise multiple connected components, each potentially containing handles  and cavities. To systematically analyze such structures, we aim to generalize the topological basis introduced in Section~\ref{sec:s4a} to 3D manifolds. This generalized basis will serve as the foundational building blocks for describing the topology of excursion sets of random fields in 3D spaces. Examples of representative topological structures that appear in this context are illustrated in Fig.~\ref{fig:3dspaces}.

{\em Topological basis}: A typical connected component, which will be an element of the topological basis in 3D, will have some numbers of handles, simply connected cavities, doubly connected cavities, and so on.
 Let $Q_{ij_0j_1...j_\ell}$ represent a connected component which has $i$ number of handles, $j_0$ number of simply connected cavities, $j_1$  number of 1-multiply connected cavities, and so on. 
Here $\ell$ denotes the largest number of multiply connectedness of the cavities. Each of the indices $j_0,j_1,j_2,...,j_\ell$ take values from 0 to $j_{\rm max}$. As in the 2D case, $j_{\rm max}$ is bounded by $q$ owing to geometrical constraints. The  value of $\ell$ is also bounded by $j_{\rm max}$. 
The set $B=\{Q_{ij_0j_1...j_\ell}\}$ 
forms a topological basis for excursion sets of 3D fields. So, a typical excursion set $\mathcal{Q}^{(\nu)}$ can be written as:
\be
\mathcal{Q}^{(\nu)} =\UU \ m_{ij_0j_1...j_\ell}(\nu) Q_{ij_0j_1...j_\ell},
\ee
where  $m_{ij_0j_1...j_\ell}(\nu)$ is the number of each $Q_{ij_0j_1...j_\ell}$. 
Note that this basis does not distinguish configurations such as links and knots. 

{\em Betti numbers}: In terms of  $m_{ij_0j_1...j_\ell}$'s, $b_k$ can be expressed as
\bea
b_0 &=& \sum_{i,j_0,j_1,...,j_\ell=0}^{j_{\rm max}} m_{ij_0j_1...j_\ell}, \label{eq:3db0}\\
b_1 &=& \sum_{i,j_0,j_1,...,j_\ell=0}^{j_{\rm max}} (i+j_1+2j_2+...+\ell j_\ell)\,m_{ij_0j_1...j_\ell},  \label{eq:3db1}\quad\ \,\\ 
b_2 &=& \sum_{i,j_0,j_1,...,j_\ell=0}^{j_{\rm max}} (j_0+j_1+...+j_\ell) \,m_{ij_0j_1...j_\ell}  \label{eq:3db2}. 
\eea
Note that $j_0$ and $i$ do not enter the prefactors of $m_{ij_0j_1...j_\ell}$ for $b_1$ and $b_2$, respectively. 
Then, $\chi$ and $\bs$ can be expressed in terms of  $m_{ij_0j_1...j_\ell}$'s as the alternating sum and sum, respectively, of the right hand sides of Eqs.~\ref{eq:3db0},~\ref{eq:3db1} and \ref{eq:3db2}.
  
Similar to the 2D case, the Betti numbers can be obtained from a generating  function $h(\vec\a)$.    
Let $\vec J \equiv (i,j_0,j_1,...,j_\ell)$ and $\vec \a \equiv (\a_{-1},\a_0,\a_1,...,\a_p)$, where $\vec\a$ is a $(\ell+2)$-tuple of real variables. Let $h(\vec\a)$ be defined as,  
\be
h(\vec\a)= \sum_{\vec J=0}^{j_{\rm max}} m_{\vec J}\,e^{-\vec{J} .\vec\a}.
\ee
Then, the Betti numbers can be expressed in terms of  $h(\vec\a)$ as,
\bea
b_0 &=&  \lim\limits_{\vec\a\to 0}\ h(\vec\a), \label{eq:3db0h}\\
b_1 &=&  \lim\limits_{\vec\a\to 0}\left[ - \frac{\partial h}{\partial \a_{-1}}    + (-1)^n \sum_{n=1}^\ell \frac{\partial^n h}{ (\partial \a_{n})^n}\right],  \label{eq:3db1h}\\
b_2 &=& - \lim\limits_{\vec\a\to 0} \ \sum_{n=0}^\ell \frac{\partial h}{\partial \a_n}  \label{eq:3db2h}.
\eea
Expressions for $\chi$ and $\bs$ can then be written down as the alternating sum and sum, respectively, of the right hand sides of Eqs.~\ref{eq:3db0h},~\ref{eq:3db1h} and \ref{eq:3db2h}.

{\em Central limit theorem for $b_0, b_1, b_2, \chi$ and $\bs$}: As in the 2D case, $b_0$, $b_1$, $b_2$, $\chi$ and $\bs$, at each $\nu$, are sums or difference of non-negative discrete random variables. They can be shown to approach Gaussian nature only in the limit of $q\to\infty$, irrespective of the nature of the field. Again $q$ can approach $\infty$ provided the power spectrum of the field does not contain physical large or small cutoff scales.

{\em $m_{ij_0j_1...j_\ell}$ as Binomial variables}: The arguments for the 2D case in section \ref{sec:s4b2} can be extended to $m_{ij_0j_1...j_\ell}$'s to model them as Binomial or Hypergeometric random variables. The statistical properties of the topological statistics  then follow as consequences of the properties of $m_{ij_0j_1...j_\ell}$'s.

When multiply connected internal cavities are less probable  we may ignore indices $j_1$ onwards, and using $j=j_0$  we simply have
\be
\mathcal{Q}^{(\nu)} =\UU \ m_{ij}(\nu) Q_{ij}.
\ee
Then we have
\bea
b_0 &=& \sum_{i,j=0}^{j_{\rm max}}  m_{ij},\\
b_1 &=&  \sum_{i,j=0}^{j_{\rm max}} im_{ij},\\ 
b_2 &=& \sum_{i,j=0}^{j_{\rm max}} jm_{ij}.
\eea

A full numerical exploration of the consequences of the topological representation is beyond the scope of this paper. We  postpone such an investigation to a future work.

\section{Summary and discussion}
\label{sec:s6}

This paper aims to broaden the understanding of topological statistics beyond their ensemble expectations. 
To this end we introduced a new representation of excursion sets of smooth random fields on 2D and 3D  manifolds 
as  unions of a topological basis. Betti numbers, Euler characteristic and the sum of Betti numbers are then expressed as summations (and differences) of the coefficients of the basis elements. This provides a connection between the random nature  of the topological statistics and the nature of the discrete random coefficients, which then sets the stage for formulating the question of their statistical properties in a well defined way. Our approach in carrying out these investigations is semi-numerical owing to our ignorance of the properties of the coefficients from first principles. 

Using the new representation, we present the conditions under which each topological statistic can be asymptotically Gaussian in the limit of 
large field of view and high resolution (keeping in mind applications to cosmological data). 
We then model the coefficients of the  topological basis elements as Binomial variables, from which the statistical natures of topological statistics can be obtained. This further enables us to examine how their character changes in different threshold regimes. 
  We conclude that they are Gaussian to good approximation, except at high positive and negative thresholds. 
The new representation of excursion sets thus serves to provide mathematical clarity on the statistical foundations of topological statistics,  thereby enhancing their interpretability  and usefulness for physical inference, particularly in cosmological applications.  
A potentially useful finding is that $\bs$ has lower statistical uncertainty compared to individual Betti numbers and Euler characteristic, as a consequence of anti-correlation between the Betti numbers. 

This work raises several open questions. A key caveat in section~\ref{sec:s4} is the assumption of statistical independence among the $m_j$'s.  The semi-analytic test carried out in section~\ref{sec:s4b4} suggests that corrections to the means and variances of $m_j$'s due to their cross-correlations are likely small. However, this observation warrants further investigation, potentially through geometrical and topological reasoning, supported by numerical validation.
Investigating the spectrum of $m_j$ versus $j$ may reveal useful physical insights, potentially relevant for cosmology. A more complete treatment of the 3D case and extensions to other Minkowski functionals - with appropriate modifications of coefficients and indices to real variables - also merit future exploration.

\begin{acknowledgments}
{The author gratefully acknowledges K. P. Yogendran and Sreedhar B. Dutta for many stimulating discussions and comments on the manuscript. Thanks are also due to Dmitri Pogosyan, Najmuddin Fakhruddin, Alok Maharana, Changbom Park, Ashoke Sen and Fazlu Rahman for helpful discussions; S. Kuriki for feedback and highlighting reference~\cite{Owada:2020}; and Nidharssan S. for the python script used in Fig.~\ref{fig:pdfs}. Lastly, the author thanks Stephen Appleby for the code for Betti number computations, as well as a careful reading of the manuscript and insightful comments.} 
\end{acknowledgments}

\def\apj{ApJ}%
\def\mnras{MNRAS}%
\def\aap{A\&A}%
\def\apjl{ApJ}
\def\aj{AJ}
\def\physrep{PhR}
\def\apjs{ApJS}
\def\jcap{JCAP}
\def\pasa{PASA}
\def\pasj{PASJ}
\def\nat{Natur}
\def\apss{Ap\&SS}
\def\araa{ARA\&A}
\def\aaps{A\&AS}
\def\ssr{Space Sci. Rev.}
\def\pasp{PASP}
\def\na{New A}
\def\prd{PRD}
\def\mathann{Math. Ann.}
\def\AAP{Adv. Appl. Prob.}

\bibliographystyle{apsrev}
\bibliography{reference}

\end{document}